\documentclass[draft]{amsart}

\theoremstyle{definition}
\newtheorem{definition}{Definition}
\numberwithin{definition}{section}
\theoremstyle{plain}
\newtheorem{theorem}[definition]{Theorem}
\newtheorem{lemma}[definition]{Lemma}

\newtheorem{corollary}[definition]{Corollary}
\theoremstyle{remark}
\newtheorem{remark}[definition]{Remark}
\numberwithin{equation}{section}
\usepackage{amssymb}

\def\halfskip{\vskip 10pt plus 1pt minus 1pt}
\def\CC{\mathbb C}
\def\CF{\mathcal F}
\def\CG{\mathcal G}
\def\DD{\mathbb D}
\def\eps{\varepsilon}
\def\id{\operatorname{id}}
\def\intt{\operatorname{int}}
\long\def\comment#1{}
\def\pr{\operatorname{pr}}
\def\too{\longrightarrow}
\def\CO{\mathcal O}
\def\Omega{\varOmega}
\def\phi{\varphi}
\def\PP{\mathbb P}
\def\QQ{\mathbb Q}
\def\wdht{\widehat}

\def\XX{\mathbb X}
\def\ZZ{\mathbb Z}

\makeatletter
\def\@makefnmark{\hbox{$\left(^{\@thefnmark}\right)$\;}}
\makeatother

\begin{document}

\title
{A remark on separate holomorphy}

\author{Marek Jarnicki}
\address{Jagiellonian University, Institute of Mathematics,
Reymonta 4, 30-059 Krak\'ow, Poland} \email{Marek.Jarnicki@im.uj.edu.pl}

\author{Peter Pflug}
\address{Carl von Ossietzky Universit\"at Oldenburg, Institut f\"ur Mathematik,
Postfach 2503, D-26111 Oldenburg, Germany}
\email{pflug@mathematik.uni-oldenburg.de}
\thanks{The research was supported by DFG grant no.~227/8-1 and
was a part of the Research Grant No. 1 PO3A 005 28, which is supported
by public means in the programme promoting science in Poland in the years 2005-2008.}



\begin{abstract} Let $X$ be a Riemann domain over $\CC^k\times\CC^\ell$. If $X$ is
domain of holomorphy with respect to a family $\CF\subset\CO(X)$, then there exists
a pluripolar set $P\subset\CC^k$ such that every slice $X_a$ of $X$ with $a\notin P$
is a domain of holomorphy with respect to the family $\{f|_{X_a}: f\in\CF\}$.
\end{abstract}

\maketitle

\section{Introduction -- Riemann domains of holomorphy}

Let $(X,p)$ be a {\it Riemann domain over $\CC^n$}, i.e.~$X$ is an $n$--dimensional
complex manifold and $p:X\too\CC^n$ is a locally biholomorphic mapping
(see~\cite{JarPfl2000} for details). We say that
two Riemann domains $(X,p)$ and $(Y,q)$ over $\CC^n$ are {\it isomorphic}
(shortly, $(X,p)\simeq(Y,q)$) if there
exists a biholomorphic mapping $\phi:X\too Y$ such that $q\circ\phi=p$. Isomorphic Riemann
domains will be always identified.

We say that an open set $U\subset X$ is {\it univalent (schlicht)} if $p|_U$ is injective.
Note that the whole domain $X$ is univalent iff $(X,p)\simeq (\Omega,\id_\Omega)$,
where $\Omega$ is an open set in $\CC^n$.

Let $f\in\CO(X)$.
For any $\alpha\in\ZZ_+^n$ ($\ZZ_+$ stands for the set of
non-negative integers) and $x_0\in X$, let $D^\alpha\!f(x_0)$ denote the
{\it $\alpha$--partial derivative of $f$ at $x_0$},
$$
D^\alpha\!f(x_0):=D^\alpha(f\circ(p|_U)^{-1})(p(x_0)),
$$
where $U$ is an open univalent neighborhood of $x_0$ and $D^\alpha$ on the right hand side
means the standard $\alpha$--partial derivative operator in $\CC^n$.
Let $T_{x_0}f$ denote the {\it Taylor series of $f$ at $x_0$},
i.e.~the formal power series
$$
\sum_{\alpha\in\ZZ_+^n}\frac1{\alpha!}D^\alpha\!f(x_0)(z-p(x_0))^\alpha,\quad z\in\CC^n.
$$

For $x_0\in X$
and $0<r\leq+\infty$ let $\PP_X(x_0,r)$ denote an open univalent neighborhood of $x_0$ such
that $p(\PP_X(x_0,r))=\PP(p(x_0),r)=$ the polydisc with center at
$p(x_0)$
and radius $r$. Let $d_X(x_0)$ denote the maximal $r$ such that $\PP_X(x_0,r)$ exists.
Put $\PP_X(x_0):=\PP_X(x_0,d_X(x_0))$.

For $f\in\CO(X)$ and $x_0\in X$, let $d(T_{x_0}f)$ denote the {\it radius of convergence
of $T_{x_0}f$}, i.e.
$$
d(T_{x_0}f):=\sup\{r>0: \text{ the series $T_{x_0}f$ is convergent in } \PP(p(x_0),r)\}.
$$
Obviously, $d(T_{x_0}f)\geq d_X(x_0)$ and $f(x)=T_{x_0}f(p(x))$, $x\in\PP_X(x_0)$.
Notice that
$$
1/d(T_{x_0}f)=\limsup_{\nu\to+\infty}\Big(\max_{\alpha\in\ZZ_+^n:\; |\alpha|=\nu}
\frac1{\alpha!}|D^\alpha\!f(x_0)|\Big)^{1/\nu}.
$$

Let $\varnothing\neq\CF\subset\CO(X)$. We say that $(X,p)$ is an {\it $\CF$--domain
of existence} if
$$
d_X(x)=\inf\{d(T_xf): f\in\CF\},\quad x\in X.
$$

We say that an $\CF$--domain of existence $(X,p)$ is an {\it $\CF$--domain of holomorphy}
if {\it $\CF$ weakly separates points in $X$}, i.e.~for any $x', x''\in X$,
with $x'\neq x''$ and $p(x')=p(x'')$, there exists an $f\in\CF$
such that $T_{x'}\!f\neq T_{x''}\!f$ (as formal power series).

From now on we assume that all considered Riemann domains are countable at infinity.

\begin{remark}[Properties of domains of holomorphy]\label{Rem1}

(a) Let $(X,p)$ be an $\CF$--domain of holomorphy and let $U\subset X$ be a
univalent domain for which there exists a domain $V\supset p(U)$ such that for
every $f\in\CF$ there exists a function $F_f\in\CO(V)$ such that
$F_f=f\circ(p|_U)^{-1}$ on $p(U)$. Then there exists a univalent domain $W\supset U$
with $p(W)=V$.

Indeed, we only need to observe that we may always assume that $(X,p)$ is realized
as a subdomain of the sheaf of $\CF$--germs of holomorphic functions
(cf.~\cite{JarPfl2000}, proof of Theorem 1.8.4) and, consequently, we may put
$$
W:=\{[(D,(F_f)_{f\in\CF})]_{\overset{z}\sim}: z\in V\}
$$
(cf.~\cite{JarPfl2000}, Example 1.6.6).

(b) (\cite{JarPfl2000}, Proposition 1.8.10)
Let $A\subset X$ be a dense subset such that $A=p^{-1}(p(A))$. Then the following
conditions are equivalent:

(i) $(X,p)$ is an $\CF$--domain of holomorphy;

(ii) $d_X(x)=\inf\{d(T_xf): f\in\CF\}$, $x\in A$, and
for any $x', x''\in A$, with $x'\neq x''$ and $p(x')=p(x'')$, there exists an $f\in\CF$
such that $T_{x'}\!f\neq T_{x''}\!f$.

(c) If $(X,p)$ is an $\CF$--domain of holomorphy, then there exists a finite or countable
subfamily $\CF_0\subset\CF$ such that $(X,p)$ is an $\CF_0$--domain of holomorphy.

Indeed, we may assume that $X$ is connected. The case where
$(X,p)\simeq(\CC^n,\id_{\CC^n})$ is trivial. Thus assume that $d_X(x)<+\infty$, $x\in X$.
Let $A\subset X$ be a countable dense subset such that $A=p^{-1}(p(A))$.
By (b), for any $x\in A$ and $r>d_X(x)$ there exists an $f_{x,r}\in\CF$ such that
$d(T_xf_{x,r})<r$, and for $x', x''\in A$, with $x'\neq x''$ and $p(x')=p(x'')$, there
exists an $f_{x',x''}\in\CF$ such that $T_{x'}\!f_{x',x''}\neq T_{x''}\!f_{x',x''}$.
Now, we may take
$$
\CF_0:=\{f_{x,r}: x\in A,\;\QQ\ni r>d_X(x)\}\cup
\{f_{x',x''}: x', x''\in A,\; x'\neq x'',\; p(x')=p(x'')\}.
$$
\end{remark}

\section{Main results -- separate holomorphy}

Let $(X,p)$ be a connected Riemann domain over $\CC^n=\CC^k\times\CC^\ell$,
$$
p=(u,v):X\too\CC^k\times\CC^\ell.
$$
Put $D:=p(X)$, $D_k:=u(X)$, $D^\ell:=v(X)$.
For $a\in D_k$ define $X_a:=u^{-1}(a)$, $p_a:=v|_{X_a}$.
Similarly, for $b\in D^\ell$, put $X^b:=v^{-1}(b)$,
$p^b:=u|_{X^b}$.

\begin{remark} For every $a\in D_k$,
$(X_a, p_a)$ is a countable at infinity Riemann domain over $\CC^\ell$.
\end{remark}

Let $\varnothing\neq\CF\subset\CO(X)$. For $a\in D_k$ define $f_a:=f|_{X_a}$,
$\CF_a:=\{f_a: f\in\CF\}\subset\CO(X_a)$,
and analogously, $f^b:=f|_{X^b}$, $\CF^b:=\{f^b:f\in\CF\}\subset\CO(X^b)$, $b\in D^\ell$.

The main result of the paper is the following

\begin{theorem}\label{Thm1}
{\rm (a)} Let $\varnothing\neq\CF\subset\CO(X)$ and assume that $(X,p)$
is an $\CF$--domain of holomorphy. Then there exists a pluripolar set $S_k\subset D_k$
such that for every $a\in D_k\setminus S_k$, $(X_a, p_a)$ is
an $\CF_a$--domain of holomorphy.

{\rm (b)} Assume that $(X,p)\simeq(D,\id_D)$, where $D\subset\CC^k\times\CC^{\ell}$ is
a fat domain {\rm(}i.e.~$D=\intt\overline D${\rm)}
and there exist sets $S_k\subset D_k$, $S^\ell\subset D^\ell$ such that:
\begin{itemize}
\item $\intt S_k=\varnothing$, $\intt S^\ell=\varnothing$,
\item for any $a\in D_k\setminus S_k$, $D_a$ is an $\CF_a$--domain of holomorphy,
\item for any $b\in D^\ell\setminus S^\ell$, $D^b$ is an $\CF^b$--domain of holomorphy.
\end{itemize}
Then $D$ is an $\CF$--domain of holomorphy.
\end{theorem}

\begin{proof} (a) By Remark \ref{Rem1}(c), we may assume that $\CF$ is finite or
countable.

Step 1. {\sl There exists a pluripolar set $P\subset D_k$ such that
for any $a\in D_k\setminus P$, $(X_a,p_a)$ is an $\CF_a$--domain of existence.}

Define $R_{f,b}(x):=d(T_xf_{u(x)})$, $f\in\CF$, $b\in D^\ell$, $x\in X^b$. Recall that
$$
1/R_{f,b}(x)=
\limsup_{\nu\to+\infty}\Big(\max_{\beta\in\ZZ_+^\ell:\; |\beta|=\nu}
\frac1{\beta!}|D^{(0,\beta)}f(x)|\Big)^{1/\nu},\quad x\in X^b.
$$
Obviously, $R_{f,b}(x)\geq d_X(x)$, $x\in X^b$.
By the Cauchy inequalities, we get
$$
\frac1{\beta!}|D^{(0,\beta)}f(x)|\leq\frac{\sup_{\PP_X(x_0,r)}|f|}{r^{|\beta|}},\quad
0<r<d_X(x_0),\;x\in\PP_X(x_0,r/2),\; \beta\in\ZZ_+^\ell.
$$
Consequently, the function $-\log(R_{f,b})_\ast$
(where ${}_\ast$ denotes the lower semicontinuous regularization on $X^b$) is
plurisubharmonic on $X^b$. Put
$$
P_{f,b}:=u(\{x\in X^b: (R_{f,b})_\ast(x)<R_{f,b}(x)\})\subset D_k.
$$
It is known that $P_{f,b}$ is pluripolar (cf.~\cite{JarPfl2000}, Theorem 2.1.41(b)).
Put
$$
R_b:=\inf_{f\in\CF}R_{f,b},\quad \wdht R_b:=\inf_{f\in\CF}(R_{f,b})_\ast.
$$
Observe that $-\log(\wdht R_b)_\ast$ is plurisubharmonic on $X^b$. Put
$$
P_b:=u(\{x\in X^b: (\wdht R_b)_\ast(x)<\wdht R_b(x)\})\subset D_k.
$$
The set $P_b$ is also pluripolar (cf.~\cite{JarPfl2000}, Theorem 2.1.41(a)).
Now let $B\subset D^\ell$ be a dense countable set. Define
$$
P:=\Big(\bigcup_{f\in\CF,\;b\in B}P_{f,b}\Big)\cup\Big(\bigcup_{b\in B}P_b\Big)\subset D_k.
$$
Then $P$ is pluripolar.

Take an $a\in D_k\setminus P$ and suppose that $X_a$ is not an $\CF_a$--domain of existence.
Then there exist a point $x_0\in X_a$ and a number $r>d_{X_a}(x_0)$ such that
$b:=v(x_0)\in B$ and $R_b(x_0)>r$. Since $a\notin P$, we have
$$
(\wdht R_b)_\ast(x_0)=\wdht R_b(x_0)=\inf_{f\in\CF}(R_{f,b})_\ast=\inf_{f\in\CF}R_{f,b}
=R_b(x_0)>r.
$$
In particular, there exists $0<\eps<d_X(x_0)$ such that
$(\wdht R_b)_\ast(x)>r$, $x\in\PP_{X^b}(x_0)$.
Since,
$$
R_b(x)=\inf_{f\in\CF}R_{f,b}(x)\geq\inf_{f\in\CF}(R_{f,b})_\ast(x)=\wdht R_b(x)\geq
(\wdht R_b)_\ast(x),
$$
we conclude that $R_b(x)>r$, $x\in\PP_{X^b}(x_0)$. Put $U:=\PP_X(x_0,\eps)$.
Hence, by the classical Hartogs lemma,
for every $f\in\CF$, the function $f\circ(p|_U)^{-1}$ extends holomorphically
to $V:=\PP(a,\eps)\times\PP(b,r)$.
Since $(X,p)$ is an $\CF$--domain of holomorphy,
by Remark \ref{Rem1}(a), there exists a univalent domain $W\subset X$, $U\subset W$, such that
$p(W)=V$. In particular, $d_{X_a}(x_0)\geq r$; contradiction.

\halfskip

Step 2. {\sl There exists a pluripolar set $P\subset D_k$ such that
for any $a\in D_k\setminus P$ the family $\CF_a$ weakly separates points in $X_a$.}

Take $a\in D_k$, $x', x''\in X_a$ with  $x'\neq x''$ and
$p_a(x')=p_a(x'')=:b$. Since $\CF$ weakly separates points in $X$,
there exists an $f\in\CF$ such that
$T_{x'}f\neq T_{x''}f$. Put $r:=\min\{d(T_{x'}\!f),\;d(T_{x''}\!f)\}$ and let
$$
P_{a,x',x''}:=\bigcap_{w\in\PP(b,r)}\{z\in\PP(a,r): T_{x'}\!f(z,w)=T_{x''}\!f(z,w)\}.
$$
Then $P_{a,x',x''}\varsubsetneq\PP(a,r)$ is an analytic subset.
For any $z\in\PP(a,r)\setminus P_{a,x',x''}$ we have
$T_{x'}f\!(z,\cdot)\not\equiv T_{x''}f\!(z,\cdot)$ on $\PP(b,r)$.

Take a countable dense set $A\subset D_k$.
For any $a\in A$ let $B_a\subset X_a$ be a countable dense subset such that
$p_a^{-1}(p_a(B_a))=B_a$. Then
$$
P:=\bigcup_{\substack{a\in A,\; x',x''\in B_a\\ x'\neq x'',\; p_a(x')=p_a(x'')}}
P_{a,x',x''}
$$
is a pluripolar set.

Fix $a_0\in D_k\setminus P$, $x'_0, x''_0\in X_{a_0}$,  with $x'_0\neq x''_0$ and
$p_{a_0}(x'_0)=p_{a_0}(x''_0)=:b_0$. Put $r:=\min\{d_X(x'_0),\;d_X(x''_0)\}$.
Let $a\in A\cap\PP(a_0,r/2)$ and $x', x''\in B_a$ be such that
$x'\in\PP_X(x'_0,r/2)$, $x''\in\PP_X(x''_0,r/2)$, $p_a(x')=p_a(x'')$. Since
$a_0\notin P$, we conclude that
$T_{x'}\!f(a_0,\cdot)\not\equiv T_{x''}\!f(a_0,\cdot)$ on $\PP(b_0,r/2)$.
Consequently,
$T_{x'_0}\!f(a_0,\cdot)\not\equiv T_{x''_0}\!f(a_0,\cdot)$ on $\PP(b_0,r/2)$, which
implies that $T_{x'_0}\!f_{a_0}\neq T_{x''_0}\!f_{a_0}$.

(b) Suppose that there exist $x_0=(a_0,b_0)\in G$ and $r>d_D(x_0)=:r_0$ such that
$d(T_{x_0}f)\geq r$, $f\in\CF$. Then $d(T_{b_0}f_a)\geq r$ for any $f\in\CF$ and
$a\in\PP(a_0,r_0)$. Consequently, $(\PP(a_0,r_0)\setminus S_k)\times\PP(b_0,r)\subset D$.
Since $\intt S_k=\varnothing$ and $D$ is fat, we conclude that
$\PP(a_0,r_0)\times\PP(b_0,r)\subset D$. Now, we see that
$d(T_{a_0}f^b)\geq r$ for any $f\in\CF$ and
$b\in\PP(b_0,r)$. Consequently, $\PP(a_0,r)\times(\PP(b_0,r)\setminus S^\ell))\subset D$.
Hence $\PP(a_0,r)\times\PP(b_0,r)\subset D$; contradiction.
\end{proof}

\begin{remark}
The following natural question arises from the discussion above: is it possible to
sharpen Theorem \ref{Thm1}(a) so that the exceptional set there is even a countable union of locally analytic sets.
The following example will show that the answer is, in general,
negative.

Let $C_1\subset \DD:=$ (the unit disc) be a compact polar set which is uncountable
(take, for example, an appropriate Cantor set). Define $C:=C_1\cup
C_2$, where $C_2:=\DD\cap\QQ^2$. Then $C$ is polar
and a countable union of compact sets.

Using  Example 2 from \cite{Ter1972},
we find a function $f:\DD\times \DD\too\CC$ with the following
properties
\begin{itemize}
\item  $f(\cdot, w)\in\CO(\DD)$,\quad $w\in \DD$,
\item $f(z,\cdot)\in\CO(\DD)$,\quad $z\in C$,
\item $f$ is unbounded near some point $(z_0,0)\in \DD\times \DD$.
\end{itemize}

Using the corollary to Lemma 8 in \cite{Ter1972}, we conclude that there is a non-empty
open set $V\subset\subset \DD$ such that $f|_{\DD\times V}\in\CO(\DD\times V)$. Set
$\CF:=\{f|_{\DD\times V}, g|_{\DD\times V}\}$, where $g\in\CO(\DD\times
\DD)$ is chosen such that $\DD\times \DD$ is the existence domain of
$g$. Denote by $(D',p)$ the $\CF$-envelope of holomorphy of $\DD\times V$. Then $p(D')\subset \DD\times \DD$.
Moreover, using the fact that $C$ is dense in $\DD$ one sees that $D'$ is univalent. Indeed, let us take a sequence
$G_j=G_j'\times G_j''\subset \DD\times \DD$, $j=1,\dots,N$, of bidiscs, $G_j\cap G_{j+1}\neq\emptyset$, and functions $f_j\in\CO(G_j)$, $j=1,\dots,N$, such
that $G_1\subset \DD\times V$, $f_1=f|_{G_1}$, and $f_j|_{G_j\cap G_{j+1}}=f_{j+1}|_{G_j\cap G_{j+1}}$, $j=1,\dots,N-1$.
We claim that then $f_N=f|_{G_N}$ which implies that $D'$ is
univalent. By induction we may assume that $f_j=f|_{G_j}$ for an $j<N$. Then for any point
$a\in C\cap G_{j}'\cap G_{j+1}'$ we have two holomorphic functions $f(a,\cdot)$ and $f_{j+1}(a,\cdot)$ on $G_{j+1}''$.
Both coincide on $G_j''\cap G_{j+1}''$, and so on $G_{j+1}''$. Now fix a $b\in G_{j+1}''$. Then $f(\cdot,b)$ and
$f_{j+1}(\cdot,b)$ are holomorphic in $G_{j+1}'$ and they coincide on $C\cap G_{j}'\cap G_{j+1}'$; hence they are
equal on $G_{j+1}''$, i.e. $f|_{G_{j+1}}=f_{j+1}$.

Set $D:=p(D')$.
Then $D$ is an $\wdht\CF$-domain of holomorphy, where
$$
\wdht\CF:=\{\wdht g:=g|_D, \wdht f:=f|_D\}.
$$
Observe that for
any $a\in C$, the functions $\wdht f(a,\cdot), \wdht g(a,\cdot)$ extend to the whole of $\DD$.

Fix $R'<R\in (0,1)$ such that $V\subset\subset \PP(0,R')$.
Suppose that there is an $a_0\in C$ with $\{a_0\}\times
\PP(0,R)\subset D$. Then there is a small open neighborhood $U\subset \DD$ of $a_0$ such that $U\times\PP(0,R')\subset D$.
In virtue of the Hartogs lemma we conclude that $f$ is holomorphic on $\DD\times\PP(0,R')$, in particular
an holomorphic extension of $f|_{\DD\times V}$, and therefore
bounded near $(z_0,0)$; a contradiction. Thus, the singular set $S_1$ for $D$ must
contain $C$.
\end{remark}

\begin{remark} Observe that Theorem \ref{Thm1}(b) need not be true if $D$ is not fat.
For example, let $D:=\DD^2\setminus\{(0,0)\}\subset\CC^2$, $\CF:=\CO(D)$.
By the Hartogs extension theorem, any function from $\CF$ extends holomorphically
to $\DD^2$. Thus $D$ is not an $\CF$--domain of holomorphy.
Observe that for any $a\in\DD\setminus\{0\}$, $D_a=\DD$ and $\CF_a=\CO(\DD)$.
Similarly, for any $b\in\DD\setminus\{0\}$, $D^b=\DD$ and $\CF^b=\CO(\DD)$.
\end{remark}

\section{Applications -- separately holomorphic functions}

Directly from Theorem \ref{Thm1} we get the following useful corollary.

\begin{corollary}\label{Cor1}
Let $D\subset\CC^k\times\CC^\ell$ be a domain, let $\varnothing\neq\CF\subset\CO(D)$ and let
$A\subset\pr_{\CC^k}(D)$. Assume that for any $a\in A$ we are given
a domain $G(a)\supset D_a$ in $\CC^\ell$ such that:
\begin{itemize}
\item for any $f\in\CF$, the function $f(a,\cdot)$ extends to an $\wdht f_a\in\CO(G(a))$,
\item the domain $G(a)$ is a $\{\wdht f_a: f\in\CF\}$--domain of holomorphy.
\end{itemize}
Let $(X,p)$ be the $\CF$--envelope of holomorphy of $D$.
Then there exists a pluripolar set $P\subset A$ such that for every $a\in A\setminus P$ we have
$(X_a,p_a)\simeq(G(a),\id_{G(a)})$.
\end{corollary}

Recall a version of the cross theorem for separately holomorphic functions with
pluripolar singularities (cf.~\cite{JarPfl2003}, Main Theorem).

\begin{theorem}\label{Thm2}
Let $D\subset\CC^k$, $G\subset\CC^\ell$ be domains of holomorphy and let
$A\subset D$, $B\subset G$ be locally pluriregular sets
\footnote{A non-empty set $A$ is said to be {\it locally pluriregular} if for any
$a\in A$, the set {\it $A$ is locally pluriregular at $a$}, i.e.~for any open neighborhood $U$ of $a$ we have $h^\ast_{A\cap U,U}(a)=0$,
where $h_{A\cap U,U}$ denotes the relative extremal function of $A\cap U$ in $U$.
For arbitrary set $A$ define $A^\ast:=\{a\in\overline A: A \text{ is locally pluriregular
at } a\}$. It is known that the set $Z:=A\setminus A^\ast$ is pluripolar. In particular,
if $A$ is non pluripolar, then $A\setminus Z$ is locally pluriregular.}\!. Consider the {\rm cross}
$$
X=\XX(A,B;D,G):=(A\times G)\cap(D\times B)
$$
and let
$$
\wdht X=\wdht\XX(A,B;D,G):=\{(z,w)\in D\times G: \omega_{A,D}(z)+\omega_{B,G}(w)<1\},
$$
where $\omega_{A,D}$ and $\omega_{B,G}$ are generalized relative extremal functions.
Let $M\subset X$ be a relatively closed set such that:
\begin{itemize}
\item for every $a\in A$ the fiber $M_a:=\{w\in G: (a,w)\in M\}$ is pluripolar,
\item for every $b\in B$ the fiber $M^b:=\{z\in D: (z,b)\in M\}$ is pluripolar.
\end{itemize}
Let $\CF=\CO_s(X\setminus M)$ denote the set of all functions
{\rm separately holomorphic on $X\setminus M$},
i.e.~of those functions $f:X\setminus M\too\CC$ for which:
\begin{itemize}
\item for every $a\in A$, $f(a,\cdot)\in\CO(G\setminus M_a)$,
\item for every $b\in A$, $f(\cdot,b)\in\CO(D\setminus M^b)$.
\end{itemize}
Then there exists a relatively closed pluripolar set $S\subset\wdht X$ such that:
\begin{itemize}
\item $S\cap X\subset M$,
\item for every $f\in\CF$ there exists an $\wdht f\in\CO(\wdht X\setminus S)$
with $\wdht f=f$ on $X\setminus M$,
\item $\wdht X\setminus S$ is a $\{\wdht f: f\in\CF\}$--domain of holomorphy.
\end{itemize}
\end{theorem}

In the proof presented in \cite{JarPfl2003} the assumption that $M$ is relatively closed
in $X$ played an important role. Observe that from the point of view of the formulation
of the above theorem, we only need to assume that all the fibers $M_a$ and $M^b$ are
relatively closed. Corollary \ref{Cor1} permits us to clarify this problem in certain
cases.

\begin{lemma}\label{Lem1}
Let $D\subset\CC^k$, $G_0\subset G\subset\CC^\ell$ be domains of holomorphy and let
$A\subset D$. Assume that for every $a\in A$ we are given a relatively
closed pluripolar set $M(a)\subset G$. Let $\CF$ denote the set of all functions
$f\in\CO(D\times G_0)$ such that for every $a\in A$, the function $f(a,\cdot)$ extends
to an $\wdht f_a\in\CO(G\setminus M(a))$. Assume that for every $a\in A$ the set $M(a)$
is singular with respect to the family $\{\wdht f_a: f\in\CF\}$
\footnote{Recall that a relatively closed pluripolar set $M\subset G$ is said to be
{\it singular} with respect to a family $\varnothing\neq\CG\subset\CO(G)$ if there is no
point $a\in M$ which has an open neighborhood $U\subset G$ such that every function
from $\CG$ extends holomorphically to $U$
(cf.~\cite{JarPfl2000}, \S\;3.4).}\!. Then there exists a
pluripolar set $P\subset A$ such that if we put $A_0:=A\setminus P$, then the set
$$
M(A_0):=\bigcup_{a\in A_0}\{a\}\times M(a)
$$
is relatively closed in $A_0\times G$.
\end{lemma}

\begin{proof} First observe that every function from $\CO(G)$ may be regarded as an
element of $\CF$, which implies that
for every $a\in A$ the domain $G(a):=G\setminus M(a)$ is a
$\{\wdht f_a: f\in\CF\}$--domain of holomorphy.

Let $(X,p)$ be the $\CF$--envelope of holomorphy of $D\times G_0$. Since $D$ and $G$ are
domains of holomorphy, we may assume that $p(X)\subset D\times G$.

By Corollary \ref{Cor1}, there exists a pluripolar set $P\subset A$ such that for every
$a\in A_0:=A\setminus P$ we have $(X_a,p_a)\simeq(G(a),\id_{G(a)})$.
Thus $p$ is injective on the set $B:=p^{-1}(A_0\times G)$ and
$p(B)=\bigcup_{a\in A_0}\{a\}\times G(a)=(A_0\times G)\setminus M(A_0)$.
Hence $p(B)=p(X)\cap(A_0\times G)$ and, consequently, $p(B)$ is relatively open in $A_0\times G$.
\end{proof}

Consequently, we get the following generalization of Theorem \ref{Thm2}.

\begin{theorem}
Let $D_0\subset D\subset\CC^k$, $G_0\subset G\subset\CC^\ell$ be domains of holomorphy
and let $A\subset D_0$, $B\subset G_0$ be non-pluripolar sets. Let
$M\subset X:=\XX(A,B;D,G)$ be such that:
\begin{itemize}
\item for every $a\in A$ the fiber $M_a$ is relatively closed pluripolar subset of $G$,
\item for every $b\in B$ the fiber $M^b$ is relatively closed pluripolar subset of $D$.
\end{itemize}
Let $\CF$ denote the set of all functions $f\in\CO(D_0\times G_0)$ such that:
\begin{itemize}
\item for every $a\in A$ the function $f(a,\cdot)$ extends holomorphically to
$G\setminus M_a$,
\item for every $b\in A$ the function $f(\cdot,b)$ extends
holomorphically to $D\setminus M^b$.
\end{itemize}
Then there exist:
\begin{itemize}
\item pluripolar sets $P\subset A$, $Q\subset B$ such that
the sets $A_0:=A\setminus P$ and $B_0:=B\setminus Q$ are locally pluriregular,
\item a relatively closed pluripolar set $S\subset \wdht X_0:=\wdht\XX(A_0,B_0;D,G)$
\end{itemize}
\noindent such that:
\begin{itemize}
\item $S\cap\XX(A_0,B_0;D,G)\subset M$,
\item for every $f\in\CF$ there exists an $\wdht f\in\CO(\wdht X_0\setminus S)$
with $\wdht f=f$ on $D_0\times G_0$,
\item $\wdht X_0\setminus S$ is a $\{\wdht f: f\in\CF\}$--domain of holomorphy.
\end{itemize}
\end{theorem}

\begin{remark} More general versions of the cross theorem (also for $N$--fold crosses)
 will be discussed in our forthcoming paper.
\end{remark}

{\bf Acknowledgment.}
The authors like to thank Professor W.~Zwonek for helpful remarks.

\bibliographystyle{amsplain}

\end{document}